\documentclass[11pt]{article}
\usepackage{amsfonts}
\usepackage{amssymb}
\renewcommand{\le}{\leqslant}\renewcommand{\ge}{\geqslant}
\newcommand{\n}{\noindent}\newcommand{\fb}[1]{\framebox{$#1$}}
\renewcommand{\,}{\kern2pt}\newcommand{\ds}{\displaystyle} 
\renewcommand{\c}[1]{\cite{#1}}\def\qq{\qquad}\def\q{\quad}
\def\pr{\n{\sl Proof.\q}}
\newcommand{\hb}[1]{\hbox{\rm#1}}
\def\bb{\bigbreak}\def\mb{\medbreak}\def\sb{\smallbreak}
\newcommand{\DD}{\mathbb{D}}\def\V{\DD}\newcommand{\CC}{{\cal C}}
\newcommand{\GG}{\mathbb{G}r_3(\gsc)}\renewcommand{\SS}{{\cal S}}
\newcommand{\UU}{{\cal U}}
\def\oa{\ol\a}\def\ob{\ol\b}\def\og{\ol\ga}\def\os{\ol\si}
\def\ot{\ol\tau}\def\ow{\ol\w}
\def\Im{\hb{Im}}
\def\Aut{{\cal A}}\def\Der{\mathfrak{d}}\def\Hom{\hb{Hom}}
\newcommand{\g}{\mathfrak{g}}
\newcommand{\gs}{\g^*}\newcommand{\gsc}{\g^*_c}
\newcommand{\C}{\mathbb{C}}\newcommand{\R}{\mathbb{R}} 
\renewcommand{\a}{\alpha}\renewcommand{\b}{\beta}
\newcommand{\ga}{\gamma}\newcommand{\Ga}{\Gamma} 
\renewcommand{\th}{\theta}
\def\De{\Delta}\def\fh{\mathfrak{h}} 
\def\qed{\q$\square$\hfill\medbreak}\newcommand{\w}{\omega}
\newcommand{\si}{\sigma} 
\newcommand{\la}{\lambda}\newcommand{\La}{\Lambda} 
\renewcommand{\L}{\La^{1,0}}\newcommand{\we}{\wedge}
\newcommand{\D}{\partial}\newcommand{\oD}{\ol\D}
\newcommand{\angl}{\langle}\newcommand{\angr}{\rangle}
\def\mto{\,\mapstochar\kern-1.8pt\joinrel{\hbox to15pt{\rightarrowfill}}\,}
\newcommand{\cir}{\hbox{\footnotesize$\circ$}} 
\def\E{\raise1pt\hbox{$\textstyle\bigwedge$}\kern-1pt}
\newcommand{\ol}{\overline}\newcommand{\op}{\oplus}
\newcommand{\rf}[1]{(\ref{#1})}\newcommand{\cO}{{\cal O}}
\newcommand{\be}{\vspace{-1pt}\begin{equation}} 
\newcommand{\ee}[1]{\label{#1}\end{equation}} 
\newcommand{\ba}{\begin{array}}\newcommand{\ea}{\end{array}} 
\outer\def\pro#1#2\par{\vskip5pt\bigbreak\noindent{\bf#1\enspace}{\sl#2} 
\par\vskip3pt\medskip}\newcommand{\NLA}{NLA}
\newcommand{\frt}[2]{\hbox{$\textstyle\frac{#1}{#2}$}} 
\begin{document}\parskip2pt\parindent15pt\mathsurround1pt 

\phantom.\vspace{-100pt}
\phantom.\hfill{\small\sf To appear in J.~Pure Appl.~Algebra} 
\vspace{60pt}

\centerline{\LARGE\bf Complex structures on nilpotent Lie algebras}
\vspace{15pt}

\centerline{\large S.\,M.~Salamon}\vspace{20pt}

\n{\bf Abstract.} We classify real 6-dimensional nilpotent Lie algebras for
which the corresponding Lie group has a left-invariant complex structure, and
estimate the dimensions of moduli spaces of such structures.\vspace{5pt}

\n AMS: 17B30, 32G05, 53C30.

\setcounter{enumii}0\setcounter{equation}0 
\subsection*{Introduction} 

Let $G$ be a simply-connected real nilpotent Lie group of dimension $m$. The
nilpotency is equivalent to the existence of a basis $\{e^1,\ldots,e^m\}$ of
left-invariant 1-forms on $G$ such that \be de^i\in\E^2\angl
e^1,\ldots,e^{i-1}\angr,\q 1\le i\le m,\ee{nil} where the right-hand side is
interpreted as zero for $i=1$. For any such $G$ with rational structure
constants, there is a discrete subgroup $\Ga$ such that $M=\Ga\backslash G$ is
a compact manifold \c{Ma}. In this case, the differential graded algebra
$\mathbb A$ of left-invariant forms on $G$ is isomorphic to the de\,Rham
algebra of $M$, and the Betti numbers $b_i$ of $M$ coincide with the
dimensions of the Lie algebra cohomology groups of $\g$ \c{N}. Property
\rf{nil} then implies that $\mathbb A$ provides a minimal model for $M$ in the
sense of Sullivan, and $\mathbb A$ cannot be formal unless $G$ is abelian and
$M$ is a torus \c{DGMS,BG,H,Ca}.

Our interest arises from the imposition of extra geometrical structures on
$M$. In particular, suppose that $m=2n$ and $J$ is a complex structure on $M$
associated to a left-invariant tensor on $G$. We shall work mainly in the
language of differential forms, and a starting point for an analysis of
low-dimensional examples is the observation that in these circumstances there
exists a non-zero closed $(1,0)$-form. More generally, there exists a basis
$(\w^1,\ldots,\w^n)$ of $(1,0)$-forms such that $d\w^i$ belongs to the ideal
generated by $\w^1,\ldots,\w^{i-1}$. This result is proved in \S1 by dualizing
the central descending series of $\g$, and is subsequently used to carry
through arguments by induction on dimension.

If the complex structure $J$ makes $G$ a complex Lie group, then the exterior
derivative of any left-invariant $(1,0)$-form has type $(2,0)$. On the other
hand, the complex structure $J$ is `abelian' if the derivative of any
invariant $(1,0)$-form has type $(1,1)$ \c{BDM,FD}. In both cases it follows
that \be d\w^i\in\E^2\angl\w^1,\ldots,\w^{i-1},\ow^1,\ldots,\ow^{i-1} \angr,\q
1\le i\le n.\ee{cn} The properties of complex structures on nilmanifolds with
a basis of $(1,0)$-forms satisfying \rf{cn} is the subject of papers by
Cordero, Fern\'andez, Gray and Ugarte \c{CFG, CFGU}. Our approach had its
origin in the realization that not all complex structures have a basis of this
form, a fact that makes the general theory all the more richer.

In \S3 we embark on a determination of nilpotent Lie algebras of dimension 6
giving rise to complex nilmanifolds, without quoting general classification
results. For this reason, our results are likely to serve for a fuller
understanding of complex structures on 8-dimensional nilmanifolds. The failure
of \rf{cn} for $i=2$ is illustrated by the assertion that there exists a
unique 6-dimensional nilpotent Lie algebra (\NLA) $\g$ for which $G$ has a
left-invariant complex structure $J$ and $b_1=2$. The cases $b_1=3,4$ are more
complicated and, in order to streamline the presentation, we precede the
classification by a number of preliminary results in \S2. Conformal structures
are used to construct canonical bases for subspaces of 2-forms in 4
dimensions, and this study leads to a dichotomy in the choice of a canonical
basis (Theorem~2.5).

A general classification of nilpotent Lie algebras exists in dimension 7 and
less, though 6 is the highest dimension in which there do not exist continuous
families. According to published classifications, there are 34 isomorphism
classes of \NLA s over $\R$, of which 10 are reducible \c{M,GK}. It is also
known which of the corresponding algebras admit symplectic structures, and in
\S4 we parametrize complex and symplectic structures on 6-dimensional \NLA s,
and discuss their deformation. This part develops the approach of \c{AGS},
though here we do not impose compatibility with a metric. Proposition~4.2
concerns the infinitesimal theory, whose non-trivial nature is illustrated by
an example of an obstructed cocycle on the Iwasawa manifold. The variety
$\CC(\g)$ of complex structures on a Lie algebra $\g$ can be regarded as the
fibre of a reduced twistor space over the corresponding Lie group, in the
spirit of \cite{BR}.

Finally, we insert our results into the classification in an appendix that was
influenced by \c{CFGU}, and completes the picture given there. Our
integrability equations \rf{mod} can in fact be solved explicitly for the
6-dimensional \NLA s listed, and a more detailed description of the spaces
$\CC(\g)$ will appear elsewhere.

\mb\n{\bf Acknowledgments.} This paper is partly a sequel to \c{AGS} and
originated in joint discussions with E.~Abbena and S.~Garbiero, without whose
collaboration it could not have been written. The author is grateful to
M.~Vaughan-Lee for writing a program in {\sc Magma} to check cohomological
data provided in the appendix, and for useful comments from S.~Console,
I.~Dotti, A.~Fino, and G.~Ketsetzis.

\setcounter{enumii}1\subsection*{1. Invariant differential forms}

Let $\g$ be a Lie algebra of real dimension $2n$. The dual of the Lie bracket
gives a linear mapping $\gs\to\E^2\gs$ which extends to a finite-dimensional
complex \be0\to\gs\to \E^2\gs\to\E^3\gs\to\cdots\to\E^{2n}\gs\to0.\ee{cx} The
vanishing of the composition $\gs\to\E^3\gs$ corresponds to the Jacobi
identity, and, conversely, any linear mapping $d\colon\gs\to\E^2\gs$ whose
composition with the natural extension $\E^2\gs\to\E^3\gs$ is zero gives rise
to a Lie algebra. Given a Lie algebra $\g$, we denote by $b_k$ the dimension
of the $k$th cohomology space of \rf{cx}, which is isomorphic to the Lie
algebra cohomology group $H^k(\g)$ \c{CE}. This dimension equals the $k$th
Betti number of $\Ga\backslash G$ for any discrete cocompact subgroup $\Ga$,
by Nomizu's theorem \c{N}.

The descending central series of a Lie algebra $\g$ is the chain of ideals
defined inductively by $\g^0=\g$ and $\g^i=[\g^{i-1},\g]$ for $i\ge1$. By
definition, $\g$ is $s$-step nilpotent if $\g^s=0$ and $\g^{s-1}\ne0$. This
condition can easily be interpreted in terms of differential forms as follows.
Define a subspaces $\{V_i\}$ of $\gs$ inductively by setting $V_0=\{0\}$, and
\[V_i=\{\si\in\gs:d\si\in \E^2V_{i-1}\},\q i\ge1.\] Of paramount importance is
$V_1= \ker d$, and the dimension of this equals $b_1$.

\pro{Lemma 1.1} $V_i$ is the annihilator of $\g^i$. 

\pr Suppose inductively that $V_i=(\g^i)^o$; this is certainly true for $i=0$.
Then $d\si\in\E^2V_i$ if and only if $d\si$ annihilates the subspace
$\g\we\g^i$, i.e. $d\si(X,Y)=-\si[X,Y]$ vanishes for all $X\in\g$ and
$Y\in\g^i$. Equivalently, $\si\in(\g^{i+1})^o$.\qed

A left-invariant almost-complex structure on a Lie group $G$ can be identified
with a linear mapping $J\colon\g\to\g$ such that $J^2=-1$. Such a structure
determines in the usual way the subspace \[\L=\{X-iJX:X\in\gs\}\] of the
complexification $\gsc$ consisting of left-invariant $(1,0)$-forms, its
conjugate $\La^{0,1}$, and more generally subspaces $\La^{p,q}$ of $\E^{p+q}
\gsc$. The almost-complex structure $J$ is said to be {\sl integrable} if
\be[JX,JY]= [X,Y]+J[JX,Y]+J[X,JY],\ee{Nij} for all $X,Y\in\g$, and in this
case the Newlander-Nirenberg theorem implies that $(M,J)$ is a complex
manifold. We shall refer to a pair $(\g,J)$ consisting of a Lie algebra and an
integrable almost-complex structure simply as a `Lie algebra with a complex
structure'. A useful reference for the study of such objects is \c{Sn}.

The equation \rf{Nij} holds if and only if $d(\L)\subseteq\La^{2,0}
\op\La^{1,1}$, and this condition gives rise to a complex \be0\to\la^{0,0}\to
\La^{0,1}\to\La^{0,2}\to\cdots\to\La^{0,n}\to0,\ee{Dol} in which each map is
the restriction of the ordinary $\oD$ operator to a space of left-invariant
forms. It follows from \rf{cx} that the annihilator \be\g^{0,1}=(\L)^o\cong
(\La^{0,1})^*\ee{ann} has the structure of a complex Lie algebra. Only if
$J[X,Y]=[JX,Y]$ for all $X,Y\in\g$, or equivalently if $d(\L)\subseteq
\La^{2,0}$, is $\g$ the the real Lie algebra underlying $\g^{0,1}$. In this
case one can unambiguously declare that $\g$ is a `complex Lie algebra'.

Let $\pi_i$ denote the projection $(V_i)_c\to\La^{0,1}$, and consider the
complex subspace \[V_i^{1,0}=(V_i)_c\cap\L=\ker\pi_i.\] Observe that
$V_i^{1,0}\op\ol{V_i^{1,0}}$ is the complexification of the largest
$J$-invariant subspace contained in $V_i$, namely $V_i\cap JV_i$. It follows
that $\dim_\R(V_i)\ge2\dim_\C(V_i^{1,0})$. Moreover, if $\g$ is $s$-step
nilpotent then $\dim_\C(V^{1,0}_s)=n$.

The following result is an immediate consequence of Lemma~1.1. 
                                                      
\pro{Lemma 1.2} The real subspace underlying $V_i^{1,0}$ is the annihilator of 
the smallest $J$-invariant subspace containing $\g^i$, namely $\g^i(J)=\g^i+ 
J\g^i$. 

The subspace $\g^i(J)$ is in fact a subalgebra of $\g$ for every $i$. This is
because $[X,JY]$ belongs to the ideal $\g^i$ for all $X\in\g^i$, and $[JX,JY]$
belongs to $\g^i(J)$ for all $X,Y\in\g^i$ by \rf{Nij}. It follows that the
ideal in the exterior algebra generated by the real and imaginary components
of elements of $V_i^{1,0}$ is differential, i.e.\ closed under $d$.  In fact,
more is true. To explain this we denote by $I(\cal E)$ the ideal generated by
a set $\{\cal E\}$ of differential forms in the complexified exterior algebra,
and we abbreviate $\ol{\w^i}$ to $\ow^i$.

\pro{Theorem 1.3} An \NLA\ $\g$ admits a complex structure if and only if 
$\gsc$ has a basis $\{\w^1,\ldots,\w^n,$ $\ow^1,\ldots,\ow^n\}$ such that
\be d\w^{i+1}\in I(\w^1,\ldots,\w^i).\ee{ideal}\vspace{-8pt}

\pr Starting from a complex structure $J$, we construct the $\w^i$ by
successively extending a basis of $V_j^{1,0}$ to one of $V_{j+1}^{1,0}$. Given
$\w^1,\ldots,\w^i$ (or nothing if $i=0$), let $j$ be the least positive
integer (dependent on $i$) such that $V_{j+1}^{1,0}$ contains an element
$\w^{i+1}$ for which $\{\w^1,\cdots,w^{i+1}\}$ is linearly independent. By
hypothesis, $\ker\pi_j=V_j^{1,0}$ is spanned by $\{\w^1,\ldots,\w^r\}$ for
some $r\le i$.  It follows that the kernel of the linear mapping
\[\E^2\pi_j\colon \E^2(V_j)_c\to\La^{0,2}\] is a subspace of
$\angl\w^1,\ldots,\w^i\angr\we\gsc$, and this space therefore contains
$d\w^{i+1}$.

Conversely, a basis $\{\w^1,\ldots\ldots,\ow^n\}$ of $\gsc$ determines an
almost-complex structure $J$ on $\g$ by decreeing $\L$ to be the span
of the $\w^i$. The integrability of $J$ then follows from the condition
\rf{ideal}.\qed

\pro{Corollary 1.4} If $\g$ is an \NLA\ with a complex structure then
$V_1^{1,0}$ is non-zero.\vskip5pt

A complex structure on a Lie algebra $\g$ is called {\sl abelian} if \be
d(\L)\subseteq\La^{1,1};\ee{ab} this is equivalent to asserting that \rf{ann}
is an abelian Lie algebra \c{BDM}. Note that the integrability condition for
$J$ is automatically satisfied; it is the special case of \rf{Nij} for which
$[JX,JY]=[X,Y]$ for all $X,Y\in\g$.

\pro{Proposition 1.5} If $\g$ is an \NLA\ with an abelian complex structure,
then there exists a basis $\{\w^1,\ldots,\w^n\}$ of $\L$ satisfying
\rf{cn}.

\pr We mimic the proof of Theorem~1.3. Suppose that $\{\w^1,\ldots,\w^i\}$ has
been found, and let $j$ be the least positive integer such that
$V_{j+1}^{1,0}$ contains an element $\w^{i+1}$ with $\{\w^1,\ldots,\w^{i+1}\}$
linearly independent. Arguing as above, but with both projections \[\pi_j
\colon(V_j)_c\to\La^{0,1},\q\pi'_j\colon(V_j)_c\to\L,\] we obtain \[
d\w^{i+1}\in I(\w^1,\ldots,\w^i)\cap I(\ow^1,\ldots,\ow^i).\] The result
follows.\qed

\setcounter{enumii}2\subsection*{2. Choice of bases}

This section develops some algebra in dimensions 4 and 6 that will form the
basis for the classification of complex structures in \S3.

First, we summarize the essential properties of conformal structures on a real
oriented 4-dimensional vector space $\V$. Fix a non-zero element $\upsilon\in
\E^4\V$. The bilinear form $\phi$ on $\E^2\V$ defined by $\si\we\tau=\phi
(\si,\tau)\upsilon$ has signature $+++---$.  This is most easily seen by
choosing a basis $\{e^1,e^2,e^3,e^4\}$ of $\V$ with $\upsilon=e^{1234}$, and
noting that $\phi$ is diagonalized by the basis consisting of the 6 elements
\[\ba{ll}e^{12}+e^{34},\qq& e^{12}-e^{34},\\ e^{13}+e^{42},& e^{13}-e^{42},\\
e^{14}+e^{23},& e^{14}-e^{23}.\ea\] (Here and in the sequel, we adopt the
abbreviation $e^{ijk\cdots}$ for $e^i\we e^j\we e^k\we\cdots$.) This all
accords with the fact that the connected component of $O(3,3)$ is
double-covered by $SL(4,\R)$ \c{Sh}.

Let $g$ ($g$ for `metric') be an inner product on $\V$. The underlying
conformal structure is the equivalence class $[g]=\{c\,g:c>0\}$. If
$\{e^1,e^2,e^3,e^4\}$ is now an orthonormal basis for $g$, then \[\La^2_+=
\angl e^{12}+e^{34},\;e^{13}+e^{42},\;e^{14}+e^{23}\angr\] is the $+1$
eigenspace of the so-called $*$ operator, and depends only on $[g]$. Clearly,
$\phi$ is positive definite on $\La_+$, and a decomposition \be\E^2\V=\La_+
\op\La_-\ee{dec} is determined by defining $\La_-$ to be the annihilator of
$\La_+$ with respect to $\phi$. In fact, the correspondence \be [g]
\longleftrightarrow\La_+\ee{bij} is a bijection between the set of conformal
classes and the set of 3-dimensional subspaces on which $\phi$ is positive
definite.

Given \rf{dec}, any 2-form $\si$ equals $\si_++\si_-$ with $\si_\pm\in
\La^2_\pm$. The 2-form $\si$ is said to be {\sl simple} if it can be expressed
as the product $u\we v$ of 1-forms. This is the case if and only if
$\si\we\si=0$, or equivalently $|\si_+|=|\si_-|$, where the norms indicate an
inner product induced on $\E^2\V$ from one in $[g]$. On the other hand,
ignoring the conformal structure, any $\si\in\E^2\V$ such that $\si\we\si\ne0$
(a `generic' element) can be written as $e^{12}+e^{34}$ relative to a suitable
basis. The distinction between simple and generic 2-forms in 4 dimensions is
crucial, and the next two lemmas illustrate their respective properties. The
first is elementary, and its proof is omitted.

\pro{Lemma 2.1} Let $\si,\tau$ be linearly-independent simple elements of
$\E^2\V$. If $\si\we\tau=0$ then there exists a basis of $\V$ such that
$\si=e^{12}$ and $\tau=e^{13}$; otherwise there exists a basis such that
$\si=e^{12}$ and $\tau=e^{34}$.

To make \rf{bij} more explicit, fix a basis $\{e^1,e^2,e^3,e^4\}$ of $\R^4$,
let $a,b,c\in\R$, $a\ne0$, and set $\De=b^2-4ac$.  Consider the 3-dimensional
subspace \be\La=\angl e^{12}+e^{34},\,e^{13}+e^{42},\,ae^{14}+be^{42}+ce^{23}
\angr\ee{La} of $\E^2\R^4$. Relative to the given basis of $\La$, the matrix of
$\phi$ is \[\left(\ba{ccc}2&0&0\\0&2&b\\0&b&2ac\ea\right),\] so
$\phi|\La$ is positive definite if and only if $\De<0$.

\pro{Lemma 2.2} Given \rf{La}, there exists a linear transformation
$\th_-,\,\th_0$ or $\th_+$ of $\R^4$ (which one depends on $\De$) such that
\[\ba{ll} \La=\th_-(\angl
e^{12}+e^{34},\,e^{13}+e^{42},\,e^{14}+e^{23}\angr),& \hb{if $\De<0$},\\[3pt]
\La=\th_0(\angl e^{12}+e^{34},\,e^{13},\,e^{14}+e^{23}\angr),& \hb{if
$\De=0$},\\[3pt] \La=\th_+(\angl e^{12}+e^{34},\,e^{13},\,e^{24}\angr),&
\hb{if $\De>0$}.\ea\]

\pr For convenience, set \[\si=e^{13}+e^{42},\q \tau=ae^{14}+be^{42}+ce^{23}.
\] We shall define the required transformations by setting $f^i=\th_{\bullet}
(e^i)$.

Suppose first that $\De<0$, so that it is possible to solve the equation
\[ax^2+bx+c=(Ax-B)^2+(Cx-D)^2\] over $\R$ with $AD-BC=1$. Define \[\ba{ll}
f^1=Ae^1+Be^2,&\q f^4=Ae^4+Be^3,\\[3pt] f^2=Ce^1+De^2,&\q f^3=Ce^4+De^3.
\ea\] This amounts to applying a transformation in $SL(2,\R)$ simultaneously
to $\angl e^1,e^2\angr$ and $\angl e^4,e^3\angr$, and so $f^{12}=e^{12}$ and
$f^{34}=e^{34}$. Moreover, $f^{13}+f^{42}=\si$, and \[\ba{rcl}f^{14}+f^{23}
&=&(A^2+C^2)e^{14}+(AB+CD)(e^{13}-e^{42})+(B^2+D^2)e^{23}\\[3pt]&=&
\tau-\frt12b\si.\ea\] Thus $\La$ has the form stated.

Suppose that $\De>0$, $ac\ne0$, and let $s,t$ be the distinct real solutions
of the equation $ax^2+bx+c=0$ with $t\ne0$. Setting \[\ba{ll}f^1=e^1+se^2,\qq
& f^2=e^1+te^2,\\[6pt]\ds f^3=e^3+\frac1t e^4,&\ds f^4=\frac ca e^3+te^4,\ea\]
gives $f^{12}=(t-s)e^{12}$, $f^{34}=(t-s)e^{34}$, and \[f^{13}=\si+\frac
sc\tau,\qq f^{24}=\frac1a(c\si+t\tau),\] and $\angl f^{13},f^{24}\angr=\angl
\si,\tau\angr$, as required. The same conclusion can be verified if $c=0$.

Suppose that $b^2=4ac$, and let $s=-b/(2a)$. This time we define \[\ba{ll}
f^1=e^1+se^2,\qq&f^2=se^1-e^2,\\[3pt] f^3=se^3+e^4,& f^4=e^3-se^4.\ea\] Then
$f^{12}=-(s^2+1)e^{12}$, $f^{34}=-(s^2+1)e^{34}$, and \[f^{13}=\frac1{2a}
(-b\si+2\tau),\qq f^{14}+f^{23}=(s^2+1)\si.\] This completes the proof.\qed

Let $\g$ be an \NLA\ of real dimension 6. Let $\{e^1,\ldots,e^6\}$ be a basis
of $\gs$ satisfying \rf{nil}, so that there exist constants $c_{jk}^i$ such
that \be de^i=\sum_{j,k<i}c^i_{jk}e^{jk}.\ee{cijk} These equations imply in
particular that $e^1,e^2\in V_1=\ker d$. If $J$ is a complex structure on
$\g$, it follows from Corollary~1.4 that we may choose the basis in such a way
that \be e^1+ie^2\in V_1^{1,0}.\ee{e12}\sb

It turns out that in dimensions 6 and less, one may arrange that each
structure constant $c^i_{jk}$ be equal to $0,1$ or $-1$. Although we do not
need to assume this fact, it will make it easy to represent the isomorphism
class of a Lie algebra $\g$ by the equations \rf{cijk}; we write $\g$ as an
$m$-tuple $(0,0,de^3,\ldots,de^m)$ abbreviating $e^{ij}$ further to $ij$. This
notation is illustrated by

\pro{Proposition 2.3} A 4-dimensional \NLA\ admitting a complex structure is 
isomorphic to $(0,0,0,12)$ or $(0,0,0,0)$.

\pr The symbol $(0,0,0,12)$ means the Lie algebra whose dual has a basis for
which $de^i=0$ for $i=1,2,3$ and $de^4=e^{12}$, and $(0,0,0,0)$ is simply an
abelian algebra. Let $\g$ be a 4-dimensional \NLA\ with a complex
structure. We may suppose that $\L$ is spanned by a closed 1-form
$\w^1=e^1+ie^2$ together with a 1-form $\w^2=e^3+ie^4$ where $de^3=Ae^{12}$
and $de^4=Be^{12}+Ce^{13}+De^{23}$ with $A,B,C,D\in\R$. From the proof of
Theorem~1.3, \[0=\w^1\we d\w^2=(D-iC) e^{123},\] so that $D=0=C$. The result
follows by applying a linear transformation of $\angl e^3,e^4\angr$.\qed

\n It is well known that there is only one other \NLA\ in 4 dimensions, a
2-step one isomorphic to $(0,0,12,13)$, and that this does not admit a complex
structure.\bb

Returning to the case of a 6-dimensional Lie algebra with complex structure,
choose a basis of $\gs$ satisfying \rf{e12}. The annihilator $\fh=\angl e^1,
e^2\angr^o$ is then a subalgebra of $\g$, as in the proof of Lemma~1.2. Its
structure can also be described in terms of the quotient $\fh^*=\g^*/\angl
e^1,e^2\angr$ with induced operators $d$ and $J$, and any basis of $\fh^*$ has
the form $\{f^i+\angl e^1,e^2\angr:1\le i\le4\}$ with $f^i\in\gs$.  Applying
Proposition~2.3, we deduce that there exists a basis of $\gs$ satisfying
\rf{cijk} with \[de^i\in I(e^1,e^2,e^{34}),\q 1\le i\le6,\] so that $de^6$
does not involve $e^{35}$ or $e^{45}$. Stated more explicitly,

\pro{Corollary 2.4} If $\g$ is a 6-dimensional NLA with a complex structure,
then $\gs$ has a basis $\{e^1,\ldots,e^6\}$ satisfying \rf{e12} and \be\ba{l}
de^3\in\angl e^{12}\angr,\\de^4\in\angl e^{12}, e^{13},e^{23}\angr,\\de^5\in
\angl e^{12},e^{13},e^{14},e^{23},e^{24},e^{34} \angr,\\de^6 \in\angl e^{12},
e^{13},e^{14},e^{15},e^{23},e^{24},e^{25},e^{34} \angr.\ea\ee{span}\bb

We shall repeatedly choose a basis satisfying both \rf{e12} and \rf{span}, and
then exploit the freedom of choice. In particular, we are at liberty to apply
a conformal transformation \be\ba{rcl} e^1&\mto&A(\ \cos t\,e^1+\sin t\,e^2),
\\[2pt] e^2&\mto&A(-\sin t\,e^1+\cos t\,e^2)\ea\ee{conf} with $A\in\R$, and
real `triangular transformations' of the form \be e^i\mto\sum_{j=1}^iA^i_j
e^j,\q i\ge3.\ee{tri}

\pro{Theorem 2.5} Let $\g$ be a 6-dimensional \NLA\ with a complex structure.
Then there exists a basis $\{e^1,\ldots,e^6\}$ of $\gs$ satisfying \rf{span}
such that either \[\hb{(I)}\q\fb{\ba{rcl}\w^1&=&e^1+ie^2\\\w^2&=&e^3+ie^4\\
\w^3&=&e^5+ie^6\ea}\qq\hb{or}\q\ \hb{(II)}\q\fb{\ba{rcl}\w^1&=&e^1+ie^2\\
\w^2&=&e^4+ie^5\\\w^3&=&e^3+ie^6\ea}\] is a basis of $\L$ satisfying
\rf{ideal}, so that \be\ba{rcl}d\w^1&=&0,\\\w^1\we dw^2&=&0,\\\w^1\we\w^2\we
dw^3&=&0.\ea\ee{1d2}

\pr Given a basis $\{e^1,\ldots,e^6\}$ of $\gs$ satisfying \rf{e12} and
\rf{span}, set $\w^1=e^1+ie^2$ and define \[p=\dim\left(\angl
e^2,e^3,e^4\angr_c\cap\L\right),\q q=\dim\left(\angl
e^2,e^3,e^5\angr_c\cap\L\right).\] Since the complexification of a real
3-dimensional subspace cannot contain two linearly independent $(1,0)$-forms,
$(p,q)$ must equal one of $(1,0)$, $(0,1)$ or $(0,0)$.

If $p=1$, there exist real 1-forms $f^1,f^2\in\angl e^2,e^3\angr$ such that
\[\w^2=f^1+if^2+ie^4\in\L.\] Modifying the definition of $e^3,e^4$ using
\rf{tri}, we may assume that $f^1=e^3$ and $f^2=0$. There must also exist
$f^3,f^4\in\angl e^2,e^3,e^5\angr$ such that \[\w^3=f^3+if^4+ie^6\in \L,\] and
we may modify $e^5,e^6$ so that $f^3=e^5$ and $f^4=0$. This gives (I).  The
equations \rf{1d2} follow from the proof of Theorem~1.3.

In case $(p,q)=(0,1)$, we can modify the definition of $e^3,e^5$ using
\rf{tri} so that $e^3+ie^5\in\L$. Then \rf{1d2} implies that \be (e^1+ie^2)\we
de^5=0,\ee{12} so $e^1\we de^5=0=e^2\we de^5$ and $de^5$ is a multiple of
$e^{12}$. This being the case, we are at liberty to swap $e^4$ and $e^5$ so as
to preserve \rf{span}, reducing us to Case (I).

Now suppose that $p=0=q$. There must exist $f^1,f^2\in\angl e^2,e^3,e^4\angr$
such that \[\w^2=f^1+if^2+ie^5\in\L,\] and we may modify $e^4,e^5$ so that
$f^1=e^4$ and $f^2=0$. There also exist $f^3,f^4\in\angl e^2,e^3,e^4 \angr$
such that \[f^3+if^4+ie^6\in\L,\] and we may modify $e^6$ so that
$f^4=0$. Modifying $e^3$ if necessary, we may also assume that $f^3=e^3+ce^4$
for some $c\in\R$. Then \[\w^3=e^3+i(e^6-ce^5)=e^3+ce^4+ie^6 -c\w^2\in \L.\]
Modifying $e^6$, this gives (II), and \rf{1d2} follows as before.\qed

\setcounter{enumii}3\subsection*{3. A 6-dimensional classification}

It is convenient to divide the following analysis into cases according to the
value of the first Betti number $b_1$. Our strategy is based on Theorem~2.5,
and we shall always work with a basis $\{e^1,\ldots,e^6\}$ as described there,
satisfying (I) or (II). Observe that (I) implies (replacing $e^5$ by $e^4$ in
\rf{12}) that $de^4$ is a multiple of $e^{12}$. But then a non-zero linear
combination of $e^3,e^4$ belongs to $V_1$, and $b_1\ge3$.

\pro{Theorem 3.1} Any 6-dimensional \NLA\ with $b_1=2$ and admitting a complex 
structure is isomorphic to $(0,0,12,13,23,14+25)$.

\pr From the last remark, we can assume that we are in Case (II). Applying a
transformation \rf{conf} and then rescaling $\w^2=e^4+ie^5$, we may suppose
that $de^3=e^{12}$ and $de^4=Ae^{12}+e^{13}$, where capital letters will
always denote real coefficients. Using the equation $d(de^5)=0$, we may write
\[de^5=Be^{12}+Ce^{13}+De^{14}+Ee^{23}.\] From \rf{1d2}, \[(C+iE-i)
e^{123}+De^{124}=0,\] so $C=D=0$ and $E=1$.  From \rf{span} and the equation
$d(de^6)=0$, we obtain \[de^6=Fe^{12}+Ge^{13}+He^{14}+K(e^{15}+ e^{24})+
Le^{23}+Me^{25}.\] The third equation of \rf{1d2} becomes \[(G+Li)(e^{1234}+
ie^{1235})-(2K+Mi-Hi)e^{1245}=0,\] so $G=K=L=0$ and $M=H$. To summarize,
\[\ba{l} de^5=Be^{12}+e^{23},\\[3pt] de^6=Fe^{12}+ H(e^{14}+e^{25}).\ea\] Now
$H\ne0$, for otherwise $de^3,de^6$ are linearly dependent and $b_1\ge3$.
Subtracting multiples of $e^3$ from $e^5$ and $e^6$, and finally rescaling
$e^6$ completes the proof.\qed

This illustrates the technique we adopt throughout this section. We first
apply basis changes that preserve the equations characterizing a hypothetical
complex structure, in order to exploit \rf{1d2}. If these equations are
verified in a particular case then a complex structure $J$ does exist. Without
further reference to $J$, we then apply linear transformations so as to
simplify the description of the real \NLA. In the hardest cases ($b_1=3,\,4$)
it is convenient to carry out this simplification {\sl before} using the full
force of \rf{1d2}, so it remains to check which of the resulting \NLA s do in
fact carry a complex structure.

\pro{Theorem 3.2} A 6-dimensional \NLA\ with $b_1=3$ admitting a complex 
structure is isomorphic to one of \[\ba{l}(0,0,0,12,13,14),\\(0,0,0,12,13,23),
\\ (0,0,0,12,14,24),\\ (0,0,0,12,13,24),\\(0,0,0,12,13+14,24),\\
(0,0,0,12,13,14+23),\\(0,0,0,12,14,13+42),\\
(0,0,0,12,13+42,14+23),\\(0,0,0,12,23,14-35).\ea\]

\pr In Case (II), applying \rf{conf}, we may assume that
$de^4=Ae^{12}+Be^{13}$. Then \rf{1d2} implies that $de^5=Ce^{12}+Be^{23}$.
Now, $B\ne0$ for otherwise $b_1\ge4$, and rescaling $\w^2=e^4+ie^5$ we may
assume that $B=1$. Let
\[de^6=De^{12}+Ee^{13}+Fe^{14}+Ge^{15}+He^{23}+Ke^{24}+Le^{25}+Me^{34}.\] The
third equation in \rf{1d2} becomes
\[(E+iH)(e^{1235}-ie^{1234})+(F-L+iK+iG)e^{1245} +M(-e^{2345}+ie^{1345})=0,\]
and so $F=L$, $K=-G$ and $E=H=M=0$. Next, $d(de^6)=0$ implies that $G=0$,
whence $de^6=De^{12}+F(e^{14}+e^{25})$. The term $De^{12}$ may be absorbed
into $Fe^{14}$ by modifying $e^4$. Replacing $-Ce^1+Ae^2+e^3$ by a `new' $e^3$
gives \[\g\cong(0,0,0,13,23,14+25),\] which is equivalent to the last one
listed in the theorem.

In Case (I), rescaling $\w^2=e^3+ie^4$, we may assume that \be de^3=0,\q
de^4=e^{12}\ee{de4} (cf.\ \rf{12}). From \rf{1d2}, $de^6$ has no term in
$e^{15}$ or $e^{25}$, so $de^5,de^6\in\E^2\DD$, where \be\DD=\angl
e^1,e^2,e^3,e^4\angr.  \ee{DD} Furthermore, $d^2=0$ implies that neither
$de^5$ nor $de^6$ has a term in $e^{34}$. Consider the following three cases,
listed in decreasing generality:\par \n(i) at least one of $de^5\we de^5$,
$de^6\we de^6$ is non-zero;\par \n(ii) $de^5\we de^5=0=de^6\we de^6$ and
$de^5\we de^6\ne0$;\par \n(iii) $de^5\we de^5$, $de^6\we de^6$ and $de^5\we
de^6$ are all zero.

In (i), ignoring the complex structure $J$ (and swapping $e^5,e^6$ if
necessary), there exists a transformation of $\DD$ of type \rf{tri} such that
$de^5=e^{13}+e^{42}$. Thus, \[de^6=Ae^{12}+Be^{13}+Ce^{14}+De^{23}+Ee^{24}, \]
and if $A\ne0$ we can eliminate $Ae^{12}$ by replacing $-Ae^1+De^3$ by a `new'
$e^3$. Subtracting a multiple of $e^5$ from $e^6$ yields \[de^6=C'e^{14}
-E'e^{42}+D'e^{23}.\] We may now apply Lemma~2.2 and its proof to find a
linear transformation of $\DD$ preserving $\angl e^1,e^2\angr$ so as to
conclude that $\g$ is isomorphic to one of \[\ba{l}(0,0,0,12,13,14+23),\\
(0,0,0,12,14,13+42),\\(0,0,0,12,13+42,14+23).\ea\] The first of these Lie
algebras is characterized by the fact that $de^5\in \E^2V_1$. Each of them
admits a complex structure with \be\L=\angl e^1+ae^2,\,e^3+be^4,\,e^5+ce^6
\angr,\ee{L} where $c=i$ and $a,b\in\C$ solve the equation \be(e^1+ae^2)\we
(e^3+be^4)\we(de^5+ c\,de^6)=0,\ee{LL} by analogy to \rf{1d2}.

We treat (ii) and (iii) by applying Lemma~2.1 to the simple wedge products
$de^4,de^5,de^6$. Since $d^2=0$ and \rf{de4} holds, $e^{34}$ cannot belong to
$d(\gs)$. Given (iii), it follows that one may apply a linear transformation
of $\DD$ preserving the flag $\angl e^1,e^2\angr\subset V_1$, and a
transformation of $\angl e^5,e^6\angr$, in order that one of the following
holds: \[\ba{l} de^5=e^{13},\q de^6=e^{14};\\[5pt]\hspace{35pt} de^5=e^{13},
\q de^6=e^{23};\\[5pt]\hspace{70pt} de^5=e^{14},\q de^6=e^{24}.\ea\] Given
(ii), $\g$ must be isomorphic to one of the algebras
\[\ba{l}(0,0,0,12,13,24),\\(0,0,0,12,13+14,24).\ea\] All five of these Lie
algebras admit complex structures of the form \rf{L} with $a=i$ and $b,c\in\C$
solving \rf{LL}.\qed

\pro{Theorem 3.3} A 6-dimensional \NLA\ with $b_1=4$ admitting a complex
structure is isomorphic to one of \[\ba{l} (0,0,0,0,12,14+25),\\
(0,0,0,0,12,13),\\ (0,0,0,0,13+42,14+23),\\ (0,0,0,0,12,14+23),\\
(0,0,0,0,12,34).\ea\] If $J$ is a complex structure on any of the last four
algebras then $V_1$ is necessarily $J$-invariant.

\pr In Case (II), we may assume that the basis $\{e^1,\ldots,e^6\}$
constructed in the proof of Theorem~2.5 also satisfies $V_1=\angl
e^1,e^2,e^3,e^4\angr$. Using \rf{12} and rescaling $\w^2$, we may assume that
$de^5=e^{12}$. Now let \[de^6=Qe^{13}+Re^{14}+Se^{15}+Te^{23}+Ue^{24}+
Ve^{25}+We^{34}.\] By applying \rf{conf} with $A=1$ and rescaling $e^4$, may
suppose that $U=0$. The third equation in \rf{1d2} becomes \[-(Q+iT)(e^{1234}
+ie^{1235})+(S-iR+iV)e^{1245}+ W(e^{1345}+ie^{2345})=0,\] whence $Q=S=T=W=0$,
$R=V$ and we take $de^6=e^{14}+e^{25}$, as required.

In Case (I), we may suppose that \rf{de4} holds. Applying \rf{conf}, we may
also suppose that $de^5$ has no $e^{23}$-component, so that
\[de^5=Ae^{12}+Be^{13}+Ce^{14}+De^{24}+Ee^{34}.\] Then \rf{1d2} implies that
$de^6\in\E^2\DD$ (see \rf{DD}), and we argue as in the proof of Theorem~3.2,
by applying transformations of $\DD$. In case (i), we may take \[\ba{l}
de^5=e^{13}+e^{42},\\[3pt] de^6=C'e^{14}-E'e^{42}+D'e^{23}.\ea\] The condition
$b_1=4$ forces $de^4=0$ so that $V_1=\DD$ is $J$-invariant, and $\g$ is
isomorphic to one of \[\ba{l}(0,0,0,0,13,14+23),\\(0,0,0,0,13+42,14+23).\ea\]
Cases (ii) and (iii) are similar, with $\g$ is isomorphic to one of
\[\ba{l}(0,0,0,0,12,34),\\ (0,0,0,0,12,13).\ea\] All four of these \NLA s
admit complex structures of the type described in \c{AGS} with
$e^5+ie^6\in\L$.\qed

\pro{Proposition 3.4} A 6-dimensional \NLA\ with $b_1\ge5$ admitting a complex 
structure is isomorphic to any one of \[\ba{l} (0,0,0,0,0,0),\\
(0,0,0,0,0,12),\\ (0,0,0,0,0,12+34).\ea\]

\pr If $\g$ is non-abelian, then $b_1=5$. Choose a basis $\{e^1,\ldots,e^6\}$
of $\gs$ such that $V_1=\angl e^1,\ldots, e^5\angr$ and \[\w^1=e^1+ie^2,\
\w^2=e^3+ie^4\in \L.\] Since $de^6$ is a real 2-form it must have type
$(1,1)$, so that $J$ is abelian and \[de^6\in\angl
e^{12},e^{34},e^{13}-e^{42},e^{14}-e^{23}\angr.\] After a change of basis, we
may arrange that $de^6=e^{12}$ or $de^6=e^{12}+e^{34}$.\qed

\setcounter{enumii}4\subsection*{4. Moduli of complex and symplectic structures}

In this section we shall make some observations regarding the spaces of
left-invariant complex and (more briefly) symplectic structures structures on
a given nilpotent Lie group or algebra.

Fix a Lie algebra $\g$ of real dimension $2n$. Let $\CC=\{J\colon\g\to\g:
J^2=-1\}$ denote the set of all almost-complex structures on $\g$, and let
\[\CC(\g)=\{J\in\CC:[JX,JY]=[X,Y]+J[JX,Y]+J[X,JY]\}\] denote the set of
complex structures on $\g$. The choice of an almost-complex structure on $\g$
gives an identification \[\CC\cong\frac{GL(2n,\R)}{GL(n,\C)},\] and in this
sense $\CC$ is independent of $\g$. The space $\CC$ has two connected
components, corresponding to a choice of orientation, and changing the sign of
$J$ flips from one to the other if $n$ is odd. We shall see below (Example
(1)) that $\CC(\g)$ can have more than two components.

An element $J$ of $\CC$ is specified by assigning a complex $n$-dimensional
subspace $\La$ of $\gsc$ such that $\La\cap\ol{\La} =\{0\}$. This realizes
$\CC$ as an open set of the Grassmannian $\GG$, and furnishes $\CC$ with a
natural complex structure. The subspace $\La$ corresponds to the space of
$(1,0)$-forms relative to $J$, so that $\L=\La$ and $\La^{0,1}=\ol{\La}$.

\pro{Proposition~4.1} Let $\g$ be a non-abelian Lie algebra of dimension
$m=2n$. Then $\CC(\g)=\CC$ if and only $\g$ has a basis
$\{e_0,e_1,\ldots,e_{m-1}\}$ such that $[e_0,e_i]=e_i=-[e_i,e_0]$ for all
$i\ge1$ (with all other brackets zero).

\pr The existence of a basis as stated is equivalent to the existence of
$\a\in\gs$ such that $d\si=\si\we\a$ for all $\si\in\gsc$, and this clearly
implies that $\CC(\g)=\CC$. Suppose conversely that $\CC(\g)=\CC$. Fix
$\si\in\gsc$ with $\si\we\os\ne0$, and extend it to a basis
$\{\si,\si_2,\ldots,\si_n\}$ for a subspace $\La$ with
$\La\cap\ol\La=\{0\}$. The integrability condition
\[d\si\we\si\we(\si_2\we\cdots\we\si_n)=0,\] valid for all extensions, implies
that $d\si\we\si\we\eta=0$ for all $\eta\in\E^{n-1}\gs$. It follows that
$d\si\we\si=0$ and $d\si=\si\we\a$ for some $\a\in\gsc$ that depends on
$\si$. This is also valid for a real $\si\in\gs$, and working with a basis
$\{e^0,\ldots,e^{m-1}\}$ of $\gs$, it is easy to see that it possible to
choose $\a\in\gs$ such that $de^i=e^i\we\a$ for all $i$.\qed

The real tangent space $T_J\CC$ can be identified with the set of endomorphims
of $\g$ that anti-commute with $J$. Upon complexification, such a mapping
reverses the type of 1-forms, and one may identify the holomorphic tangent
space $T_J^{1,0}\CC$ with $\Hom(\L,\La^{0,1})$. To make this more explicit,
fix a basis $\{\w^1,\ldots,\w^n\}$ of $\L$. If $\tau^i(t)$ is a path in $\L$
with $\tau^i(0)=0$ and $\dot\tau^i(0)=\si^i$ for each $i$, then
\be\La_t=\angl\w^1+\ot^1(t),\ldots,\w^n+\ot^n(t)\angr\ee{Lat} is (for
sufficiently small $t$) a path in $\GG$ with tangent vector
\be\dot\th\colon\;\w^i\mto\os^i\ee{wos} at $\La_0$.

Given a complex structure $J$, \rf{Dol} gives rise to another complex \be0\to
(\La^{0,1})^*\to\Hom(\L,\La^{0,1})\stackrel{\oD}\to\Hom(\L,\La^{0,2})\to
\cdots\to\Hom(\L,\La^{0,n})\to0.\ee{tDol} Let $K$ denote the $J$-invariant
subspace of $T_J\CC$ underlying $\ker\oD$.

\pro{Proposition~4.2} If $J$ is a smooth point of $\CC(\g)$, then the tangent
space $T_J\CC(\g)$ to $\CC(\g)$ is a subspace of $K$.

\pr Suppose that \[\oD\w^i=\sum_{j=1}^n\w^j\we\oa_j^i,\qq\oa_j^i\in\La^{0,1}.
\] The almost-complex structure defined by \rf{Lat} is integrable if and only
if $d(\w^i+\ot^i(t))$ has no $(0,2)$-component for all $i$, or equivalently,
\be d(\w^i+\ot^i(t))\we (\w^1+\ot^1(t))\we\cdots\we(\w^n+\ot^n(t))=0,\q
i=1,\ldots,n.\ee{mod} If $\dot\th$ is tangent to $\CC(\g)$, this must hold to
first order in $t$, and so
\[(\oD\os^i-\sum_{j=1}^n\os^j\we\oa^i_j)\we\w^{12\cdots n}=0,\q
i=1,\ldots,n.\] This equation can be expressed as \be\oD\os^i-\dot\th(\oD\w^i)
=0,\ee{tang} which is equivalent to asserting that $\oD\dot\th=0$ in
\rf{tDol}.\qed
            
Now let $\g$ be nilpotent, and suppose that $M=G/\Ga$ is an associated
nilmanifold. The sequence \rf{tDol} is a subcomplex of the ordinary Dolbeault
complex of $M$ tensored with the holomorphic tangent bundle
$T=T^{1,0}(M,J)$. The quotient of $\ker\oD\subseteq\Hom(\L,\La^{0,1})$ by
$\oD((\L)^*)$ can be identified with the subspace of invariant classes in the
sheaf cohomology space $H^1(M,\cO(T))$. Theorem~1.3 provides a holomorphic
section $\w^{12\cdots n}$ of the canonical bundle $\E^nT^*$, and Serre duality
implies that \[ H^p(M,\cO(T))\cong H^{n-p}(M,\cO(T^*))^*=H^{1,n-p}(M)^*.\] It
follows that the complex dimension of $\CC(\g)$ does not exceed \be
n-h^{1,n}(\g)+h^{1,n-1}(\g)=n-h^{n-1,0}(\g)+h^{n-1,1}(\g)\ee{hod} in which the
Hodge numbers are computed by tensoring the finite-dimensional complex
\rf{Dol} with $\La^{p,0}$. For the deformation problem it is therefore
important to know under what circumstances the terms of \rf{hod} coincide with
the ordinary Hodge numbers of $M$, computed with forms that are not
necessarily invariant.

The equation \rf{tang} can be generalized by examining the component of the
left-hand side of \rf{mod} of type $(n,2)$. The latter also involves
second-order terms and yields the integrability equation \[\oD\th=
\frt12[\th,\th]\] for the homomorphism $\th\colon\;\w^i\mto\ot^i(t)$. Both
sides of the equation belong to $\Hom(\L,\La^{0,2})$ and the bracket is
defined relative to holomorphic sections of $T=(\L)^*$. A necessary condition
for $\dot\th\in\ker\oD$ to represent a tangent vector to $\CC(\g)$ is
therefore that $[\dot\th,\dot\th]$ be zero in cohomology (see
\cite[Theorem~5.1 and (5.86)]{Ko}). Lemma~4.3 below exhibits a situation in
which the vanishing of this primary obstruction is not automatic, so that
$\dim_\C\CC(\g)$ is strictly less than \rf{hod}.

Let $\Aut(\g)$ denote the group of automorphisms of the Lie algebra $\g$, so
that an element of $\Aut(\g)$ is a bijective linear mapping $f\colon\g\to\g$
satisfying $f([X,Y])=[f(X),f(Y)]$. Passing to the dual $\gs$, such an $f$
induces a chain mapping of the complex \rf{cx} relative to the usual action of
the general linear group on forms, and preserves the filtration $(V_i)$. If
\rf{Nij} holds and $f\in\Aut(\g)$ then $f\cir J\cir f^{-1}$ satisfies \rf{Nij}
in place of $J$, and in this way we obtain an action of $\Aut(\g)$ on
$\CC(\g)$. The tangent space to the orbit is determined by applying the Lie
algebra $\Der(\g)$ of $\Aut(\g)$ to $J$. An element $f$ of $\Der(\g)$ induces
a chain mapping of \rf{cx}, relative now to its action as a derivation on
exterior forms. For example, if $d\a=\b\we \ga$ and $f\in\Der(\g)$, then
$d(f\a)=f\b\we\ga+\b\we f\ga$. In the notation of \rf{wos}, the element of
$T_J\CC$ determined by $f$ is given by $\os^i=(f\w^i)^{0,1}$.\mb

Now suppose that $n=3$. In accordance with Theorem~1.3, we may write
\[\ba{l}\oD\w^1=0,\\[2pt]\oD\w^2=\w^1\we\oa,
\\[2pt]\oD\w^3=\w^1\we\ob+\w^2\we\og,\ea\] where $\a,\b,\ga\in\La^{1,0}$,
and \[\D\a=0,\q\D\b=\a\we\ga,\q\D\ga=0.\] Then \rf{tang} reduces to
\be\fb{\ba{l}\D\si^1=0\\[3pt]\D\si^2=\si^1\we\a\\[3pt]
\D\si^3=\si^1\we\b+\si^2\we\ga\ea}\ee{dos}

\n The solution space to these equations has dimension \rf{hod}, and this
provides an upper bound for the dimension of $\CC(\g)$. An obvious solution is
obtained by taking $\si^1=\si^2=0$ and $\os^3\in\ker\oD$. Indeed, if $\si$ is
a {\sl closed} $(1,0)$-form, then the complex structure associated to
\[\La_t=\angl\w^1,\w^2,\w^3+t\os \angr,\q t\in\C,\] belongs to the $\Aut(\g)$
orbit of $J$. In any case, $\CC(\g)$ never contains isolated points. 

\mb\n{\sl Examples~(1)}\q Let $\g$ denote the real Lie algebra
$(0,0,0,0,13+42,14+23)$, and define $\w^1,\w^2,\w^3$ by Theorem~2.5(I). Then
\[\left\{\ba{l}d\w^1=0,\\d\w^2=0,\\d\w^3=\w^{12}\ea\right.\] and both
\be\ba{l} J_0\colon\q\L=\angl\w^1,\w^2,\w^3\angr,\\[3pt]
         J_1\colon\q\L=\angl\w^1,\ow^2,\ow^3\angr\ea\ee{J01} are complex
structures on $\g$. Whilst $(\g,J_0)$ is the complex Heisenberg algebra, the
integrability of $J_1$ follows from the fact that
$d\w^3\in\La^{1,1}$. Theorem~3.3 asserts that any complex structure on $\g$
must preserve $V_1=\DD$, and therefore determines an orientation on this
4-dimensional space. Both $J_0,J_1$ induce the same overall orientation on
$\g$ but different ones on $V_1$, so $J_0,J_1,-J_0,-J_1$ lie in distinct
connected components of $\CC(\g)$. It can be shown that there are no other
components.

Let $f\in\Aut(\g)$ and suppose that $f$ preserves the overall orientation of
$\g$. Since $\Im\,d$ is spanned by the real and imaginary components of
$\w^{12}$, the restriction of $f$ to $V_1$ commutes with $J_0$, and we may
write \[f=\left(\ba{cc} \vphantom{\int_p^l}\ B\ &\b\\\ 0\ &\det B\ea\right),\]
where $B\in GL(2,\C)$ and $\b\in\C^2$. The orbit of $J_0$ under $\Aut(\g)$ is
complex 2-dimensional, whilst it is easy to check from \rf{mod} that the
connected component of $\CC(\g)$ containing $J_0$ is an open subset of a
smooth quadric in $\C^7$. The full deformations of $J_0$ were considered by
Nakamura \cite{Na} who showed that the Hodge number $h^{1,2}$ jumps up in
value at the point $J_0$. This will lead to a non-trivial behaviour of the
Fr\"ohlicher spectral sequence over the component of the moduli space
$\CC(\g)$ containing $J_0$.

\sb\n{\sl(2)}\q The stabilizer of $J_1$ in $\Aut(\g)$ corresponds to the
subgroup $S=\C^*\times \C^*$ of diagonal matrices in $GL(2,\C)$ that fix the
individual subspaces $\angl\w^1\angr$, $\angl\w^2\angr$. The orbit of $J_1$ is
thus isomorphic to the complex 4-dimensional space $\Aut(\g)/S$, whilst the
component of $\CC(\g)$ containing $J_1$ is again of complex dimension 6. On
the other hand, $\ker\oD$ has dimension 7, a situation resolved by

\pro{Lemma 4.3} Relative to $J_1$, the element \[\dot\th\colon\q\left\{
\ba{l}\w^1\mto\ow^3\\\ow^2\mto\,0\\\ow^3\mto\,0\ea\right.\] belongs to
$\ker\oD$, but is not tangent to $\CC(\g)$.

\pr Since $J_1$ is an abelian complex structure, $\oD$ annihilates any
$(0,1)$-form, and it follows from \rf{tang} that $\oD\dot\th=0$. By above, any
complex structure on $\g$ must preserve $V_1$, so if $\dot\th$ were the
tangent vector to a genuine deformation then $\dot\th(\w^1)\in\angl\ow^1,
\w^2\angr$, which is false.\qed 

\n{\sl(3)}\q Let $\fh$ be the Lie algebra characterized by Theorem~3.1. In the
notation of Theorem~2.5(II), $\fh$ has a complex structure for which $\L$ is
spanned by \[\left\{\ba{l} d\w^1=0,\\[3pt]d\w^2=\frac12(\w^1\we\w^3+\w^1\we
\ow^3), \\[3pt] d\w^3=\frac12i(\w^1\we\ow^1+\w^1\we\ow^2-\w^2\we\ow^1).\ea
\right.\] In the above notation, $\a=\frac12\w^3$, $\b=-\frac12i(\w^1+\w^2)$
and $\ga=\frac12i\w^1$. Now, $\oD(\La^{0,1})=\angl\ow^{13}\angr$ and
$\ker\oD=\angl\ow^1,\ow^3\angr$. The first equation in \rf{dos} implies that
$\si^1=a\w^1+b\w^3$ with $a,b\in\C$, but then the last equation gives $b=0$
and $\si^2+a\w^2\in\angl\ow^1,\ow^3\angr$. The middle equation gives $a=0$, so
that $\si^1=0$, which is consistent with Corollary~1.4, and the $\w^2$
component of $\si^3$ is constrained. Hence, $T^{1,0}_J(\CC)\cong\C^4$, a fact
confirmed by the calculation of $h^{2,0}(\g)=1$ and $h^{2,1}(\g)=2$.\bb

The classification of 6-dimensional nilpotent Lie groups admitting
left-invariant symplectic forms has been carried out by Goze and Khakimdjanov
\c{GK} (the author found an earlier version of this work valuable for the
classification given in the appendix). Let $G$ be a Lie group of dimension
$2n$. A left-invariant symplectic form can be identified with an element
$\si\in\E^2\gs$ such that $d\si=0$ and $\si^n\ne0$. The existence question
therefore reduces to an examination of
\[L(\g)=\ker(d\colon\E^2\gs\to\E^3\gs).\]

Let $\SS$ denote the set of non-degenerate 2-forms on $\g$. This is an open
subset of $\E^2\gs$, though fixing one form gives an identification
\[\SS\cong\frac{GL(2n,\R)}{Sp(2n,\R)}.\] The set \[\SS(\g)=\SS\cap L(\g)\] of
left-invariant symplectic forms on $G$ is therefore a (possibly empty) open
subset of $L(\g)$. If $\si\in\SS(\g)$ then the tangent space $T_\si\SS(\g)$
can be identified with $L(\g)$ itself, and \be\dim\SS(\g)=b_2+\hb{rank}
(d\colon\gs\to\E^2\gs).\ee{rank}

Now suppose that $\g$ is nilpotent of real dimension 6. In practice, the
following observations help determine whether or not $\SS(\g)$ is empty. If
$\{e^1,\ldots,e^6\}$ is a basis of $\gs$ satisfying \rf{cijk}, given $\si
\in\SS(\g)$, we may write \be\si=e^6\we f^1+e^5\we f^2+\xi,\q\xi\in\E^2\DD,
\ee{si} with $\DD$ as in \rf{DD}, $f^1\in(\DD\op\angl e^5\angr$ and $f^2\in
\DD$. Hence, \[ 0=d\si=de^6\we f^1-e^6\we df^1-e^5\we df^2+\eta,\q\eta\in
\E^3\DD,\] and it follows that $df^1=0$. In this way, the choice of basis
determines a mapping $\SS(\g)\to V_1$ defined by $\si\mto f^1$.

Moreover, if \[de^6=e^5\we f^3+\eta,\q\nu\in\E^2\DD,\] then \be df^2=
-f^1\we f^3,\ee{dsi} and \[\nu\we f^1+de^5\we f^2\in d(\E^2\DD).\] These
equations provide an informal algorithm for determining symplectic forms.

\mb\n{\sl Example~(4)}\q The \NLA\ $\g\cong(0,0,12,13,14,34+52)$ does not
admit a symplectic form $\si$. For, given $\si$ as in \rf{si} with
$f^1=Ae^1+Be^2\in V_1$, \rf{dsi} yields $df^2=-Ae^{12}$ and
$f^2=Ce^1+De^2-Ae^3$. But then \[2Ae^{134} + Be^{234}-De^{124}+d\xi=0,\] which
implies that $A=B=0$, impossible. Observe that, as a consequence of
Theorem~3.1, $\g$ does not admit a complex structure either.\bb

The moduli spaces $\CC(\g)$ and $\SS(\g)$ give rise to subsets of the
homogeneous space \[\UU\cong\frac{GL(2n,\R)}{U(n)},\] that parametrizes
Hermitian structures on $\R^{2n}$. Indeed, there is a double fibration
\[\ba{ccccc}&&\UU&&\\\raise5pt\hbox{$\pi_1$}&\swarrow&&\searrow&\raise5pt
\hbox{$\pi_2$} \\\CC&&&&\SS,\ea\] and the inverse images $\pi^{-1}(\CC(\g))$,
$\pi^{-1}(\SS(\g))$ correspond to the sets of left-invariant Hermitian and
almost-K\"ahler metrics respectively on $G$. When $n=3$, the manifolds
$\CC,\UU,\SS$ have real dimension $18,27,15$ respectively, and the table below
shows that in some cases these inverse images each have relatively small
codimension in $\UU$. On the other hand, the fact that no non-abelian \NLA\
can admit a K\"ahler metric \c{BG,Ca,CFG} implies that
\[\pi_1^{-1}(\CC(\g))\cap \pi_2^{-1}(\SS(\g))=\emptyset,\] unless $\g$ is
abelian.

\setcounter{enumii}5\subsection*{5. Appendix}

Although \S2 contains some lemmas relevant to the theory of 6-dimensional \NLA
s, we have not so far resorted to their full classification. However, in this
appendix we tabulate the list given by Magnin \c{M} (based on an earlier one
of Morosov \c{Mo}), and insert the results of this paper.

The Lie algebras appear lexicographically with respect to $(b_1,b_2,6-s)$,
where $b_i=\dim H^i(\g)$ and $s$ is the step length (see \S1). It is easy to
check that each entry satisfies the Jacobi identity $d^2=0$, and the choice of
sign is often important. For example, with reference to the first,
$(0,0,12,13,14+23,34-52)$ does not describe a Lie algebra since (in informal
notation) $d(34-52)=124-142\ne0$. The Betti number $b_3$ may be computed from
$b_1,b_2$ by means of the formula $b_3=2(b_2-b_1+1)$, which expresses the
vanishing of the Euler characteristic (for example, of an associated
nilmanifold $G/\Ga$). Dixmier \cite{Di} showed that any \NLA\ has $b_i\ge2$
for all $i$.

The column headed $\op$ indicates the dimensions of irreducible subalgebras in
case $\g$ is not itself irreducible. The numbers in the sixth column refer to
the list or irreducible algebras in \c{M}, and those in the seventh to the
table of \c{CFGU}. Referring to the latter, the pairs labelled 13,\,15 and
2,\,5 define the same Lie algebra over $\C$ as do the two 19's and the two
26's.

The number in the eighth column is the upper bound for the complex dimension
$\dim_\C\CC(\g)$ determined by \rf{hod}. It is computed (na\"\i vely) by
picking a particular complex structure $J$ on $\g$ and determining the
solution space to \rf{dos}. When $\g=(0,0,0,0,13+42,14+23)$, this solution
space has dimension 6,\,7 for the respective structures \rf{J01}, though both
components of $\CC(\g)$ have dimension 6. Indeed, the condition that any
complex structure preserves $V_1$ imposes two linear conditions, and the
equations describing $\CC(\g)$ incorporate an extra quadratic constraint.

An asterisk means that $\g$ does not admit a complex structure of type (I) in
Theorem~2.5, though the one with $b_1=4$ does admit a structure satisfying
condition \rf{cn}.

The last column lists the exact real dimension $\dim_\R\SS(\g)$, which is
readily computed using \rf{rank}. The latter implies that all connected
components of $\SS(\g)$ have the same dimension, which depends only on the
structure of the Lie algebra over $\C$. In this sense, the symplectic case is
more straightforward.\bb

As a final application, we quote a result that extends Example~4.

\pro{Theorem 5.1} The 6-dimensional \NLA s admitting neither complex nor 
symplectic structures are \[\ba{l} (0,0,12,13,14+23,34+52),\\
(0,0,12,13,14,34+52),\\(0,0,0,12,13,14+35),\\(0,0,0,12,23,14+35),\\
(0,0,0,0,12,15+34).\ea\]

\n The first two algebras in this list realize the minimum values $b_i=2$ for
$1\le i\le 6$.

\vfil\eject\begin{center}
{\large\bf Six-dimensional real nilpotent Lie algebras}\vspace{35pt}

\def\no{\hb{\bf--}}\[\ba{|c|c|c|l|c|c|c|c|c|}\hline
b_1&b_2&6\!-\!s&\hb{\qq Structure}&\op&\c{M}&\c{CFGU}&
\dim_\C\CC(\g)&\dim_\R\SS(\g)\\\hline
2&2&1&(0,0,12,13,14+23,34+52)&  &22&32&\no & \no\\
2&2&1&(0,0,12,13,14,34+52)&  &21&31&\no & \no\\
2&3&1&(0,0,12,13,14,15)&  &2&28&\no & 7\\
2&3&1&(0,0,12,13,14+23,24+15)&  &20&30&\no& 7\\
2&3&1&(0,0,12,13,14,23+15)&  &19&29&\no & 7\\
2&4&2&(0,0,12,13,23,14)&  &11&23&\no & 8\\
2&4&2&(0,0,12,13,23,14-25)&  &18&26&\no& 8\\
2&4&2&(0,0,12,13,23,14+25)&  &18&26&\>4^*& 8\\\hline
3&4&2&(0,0,0,12,14-23,15+34)&  &16&27&\no& 7\\
3&5&2&(0,0,0,12,14,15+23)&  &17&25&\no& 8\\
3&5&2&(0,0,0,12,14,15+23+24)&  &15&24&\no& 8\\
3&5&2&(0,0,0,12,14,15+24)&1+5&&22&\no& 8\\
3&5&2&(0,0,0,12,14,15)&1+5&&21&\no& 8\\
3&5&3&(0,0,0,12,13,14+35)&  &13&18&\no&\no\\
3&5&3&(0,0,0,12,23,14+35)&  &14&19&\no& \no\\
3&5&3&(0,0,0,12,23,14-35)&  &14&19&\>4^*&\no\\
3&5&3&(0,0,0,12,14,24)&1+5&&16&6&\no\\
3&5&3&(0,0,0,12,13+42,14+23)&  &10&15&6& 8\\
3&5&3&(0,0,0,12,14,13+42)&  &9&14&5 & 8\\
3&5&3&(0,0,0,12,13+14,24)&  &8&13&5 &8\\
3&6&3&(0,0,0,12,13,14+23)&  &6&11&5 & 9\\
3&6&3&(0,0,0,12,13,24)&  &7&12&5 & 9\\
3&6&3&(0,0,0,12,13,14)&  &1&10&5 & 9\\ 
3&8&4&(0,0,0,12,13,23)&  &3&7&6 & 9\\\hline
4&6&3&(0,0,0,0,12,15+34)&  &12&20&\no&\no\\
4&7&3&(0,0,0,0,12,15)&1\!+\!1\!+\!4&&17&\no& 9\\
4&7&3&(0,0,0,0,12,14+25)&1+5&&9&\>5^*& 9\\
4&8&4&(0,0,0,0,13+42,14+23)&  &5&5&6& 10\\
4&8&4&(0,0,0,0,12,14+23)&  &4&4&6 & 10\\
4&8&4&(0,0,0,0,12,34)&3+3&&2&6 & 10\\
4&9&4&(0,0,0,0,12,13)&1+5&&6&6 & 11\\\hline
5&9&4&(0,0,0,0,0,12+34)&1+5&&3&6 & \no\\
5&11&4&(0,0,0,0,0,12)&1\!+\!1\!+\!1\!+\!3&  &8&7 & 12\\\hline
6&15&5&(0,0,0,0,0,0)&1+\cdots+1&  &1&9 & 15\\\hline\ea\]
\end{center}\vfil\eject

\renewcommand{\thebibliography}{\list{{\bf\arabic{enumi}.\hfil}}
{\settowidth\labelwidth{18pt}\leftmargin\labelwidth\advance
\leftmargin\labelsep\usecounter{enumi}}\def\newblock{\hskip.05em}
\sloppy \sfcode`\.=1000\relax}
\newcommand{\bi}{\vspace{-4pt}\bibitem}

\n Mathematical Institute, 24--29 St Giles, Oxford, OX1 3LB, UK
\par\n{\tt salamon@maths.ox.ac.uk}
\enddocument